%% file: 1.tex
\documentclass[13pt]{amsart}
\usepackage{amssymb,amscd,epic, eepic, newcommand}

\newtheorem{theorem}{Theorem}[section]
\newtheorem{proposition}{Proposition}[section]
\newtheorem{definition}{Definition}[section]
\newtheorem{lemma}{Lemma}[section]
\newtheorem{corollary}{Corollary}[section]

\numberwithin{equation}{section}
\numberwithin{figure}{section}

\begin{document}
\title[ Fusion product of co-adjoint orbits]{ 
Fusion product of co-adjoint orbits}
\author[S. Chang]{Sheldon X. Chang}
\thanks{Supported in part by a Sloan Fellowship}
\maketitle
\section{Introduction}
In this paper, we introduce the fusion product of a generic pair of   co-adjoint orbits. 
This  construction provides  the geometric dual object  to the 
  product in  Verlinde fusion algebra.
The latter  is a quantum deformation of the standard tensor product,
 the fusion product constructed here  is  the corresponding deformation of
 the Cartesian product.

Let $G$ be a connected and simply-connected compact  simple Lie group,
$P_+$ be the  set of dominant integral  weights, $\theta$ be the highest root.
Suppose  $\lda\in P_+$,
 and denote by $[\lda]$   the irreducible $G$-module defined by $\lda$.
Let  $(\cdot|\cdot)$ be the invariant bilinear form on $\fg^*,\fg$,
 normalized so that $(\theta|\theta)=2$.
 For $k\in \Z_+$, define
$$P^k_+=\{\lda\in P_+|\, (\theta|\lda)\le k\}.$$

Denote by $\cLG$ the central extension of the loop group $LG$.

The set  $ P^k_+$ has 1-1 correspondence with the set of highest weight
$\cLG$-modules at level $k$.
As a consequence of    conformal field theory, one obtains 
a product on the set of integrable highest weight representations at each
level.
This product implies something  new for the representations of $G$ itself.
Namely  it induces the fusion product   $[\lda]\qtensor[\ldaa]$, 
for a pair $[\lda],[\ldaa]$ whenever  $ \lda,\ldaa\in P^k_+$.

The fusion tensor has the property that 
all the dominant integral  weights appearing in it 
are in the set   $P^k_+$; and when the level is high enough, i.e., $k\ge
(\lda+\ldaa|\theta)$,  the fusion 
tensor agrees with the standard one. 
It was conjectured by Verlinde  
 the character
of $[\lda]\qtensor[\ldaa]$ satisfies the following:
 $$\chi_{[\lda]\qtensor[\ldaa]}
=\chi_{\lda}\cdot \chi_{\ldaa} \onexplatk$$
where $M^*$ is the dual of long root lattice and $h^\vv$ is the 
dual Coxeter number.
Because the left side has only $\chi_c, c\in P_+^k$ in its expansion; and
$\{\chi_c|c\in P_+^k\}$ is an orthonormal basis of  functions on the set $\explife$, with
repsect to a suitable measure,  the above uniquely determines 
the function $\chi_{[\lda]\qtensor[\ldaa]}$.

This paper will give a proof of this  identity.
The more interesting  result here is the concrete realization of  
$[\lda]\qtensor[\ldaa]$ as holomorphic sections over the fusion product.
 
To understand the fusion product in constrast with the Cartesian product,
 we first recall the role played by the latter
in the representations of $G$.

Let $M_{\lda},M_\ldaa$ be the co-adjoint
orbits passing $\lda,\ldaa$.  Let $L_\lda$ be the line bundle over $M_\lda$
defined by the character $\lda$. Correspondingly, there is $L_\ldaa$ over
$M_\ldaa$. 

It is well known, through  Borel-Weil theory,  $[\lda]$ can be realized as
 $H^0(M_\lambda, L_\lambda)$. 
Obviously then  $[\lda]\otimes [\ldaa]$ can be realized as
$H^0(M_{\lambda}\times M_{\lambda'}, L_{\lambda}\otimes L_{\lambda'} )$. 

Hence arises the  question  as how to describe the
 new represenation
$[\lda]\qtensor[\ldaa]$ in a geometric manner. 

 For a generic pair $\lda,\ldaa\in P_+^k$, the fusion
product $Y$  is obtained by performing certain surgery on the product of conjugacy classes
$$\Ad_Ge^{2\pi i \lda/k}\times \Ad_Ge^{2\pi i \ldaa/k}.$$
The  fusion product is a   holomorphic  $G$-orbifold with a orbifold line
bundle $V$, and  has the following property:
$$ \tr(s| H^0(Y,V^k))=\RR(H^0(Y,V^k))(s)=\chi_\lda(s)\cdot\chi_\ldaa(s) \onexplatk.$$

 One essential ingredient in verifying the above results 
  is the new fixed point formula for loop group action proved in [C2].
The explicit calculation of the weights at the fixed points of  the fusion
product,  which uses   properties of affine roots transformed by  affine Weyl
group, is presented here for the first time.

The first step  is to construct moduli space $\cMab$ of flat connections over three-holed
Riemann sphere, with fixed holonomies around two circles and no
constraint on the third one.  It is a
$LG$-space, and its quotient by the nilpotent subgroup $LG^{\C+}$ 
has  an $T$-orbifold as its compactification, following the main result in
[C1]. The fusion product is defined as  $Y=G\times_T X_N$, which 
can be described in terms of the product of  conjugacy classes.  
Using the fixed point formula of [C2], we calculate the equivariant
Riemann-Roch of the pair $(Y,V^k)$ 
 which agrees with what Verlinde conjectured.  

 There are much in common between $Y$ and $\Ad_G \lda\times \Ad_G \ldaa$.
The interior fixed points of $Y$ has 1-1 correspondence with the set of fixed
points  on the
Cartesian product. Although the weights are different. It turns out
the contributions   of the corresponding fixed points to Riemann-Roch differ by a factor which is 1 on the set $\explife$, as asserted by  Theorem 4.1.

Another  important difference
is that for the standard tensor, there exists one highest weight, which is
not true for the fusion tensor. The figure 4.2 illustrates this point.

For a generic pair $\lda,\ldaa$,  we remark when the level $k\ge
(\lda+\ldaa|\theta)$,   one does not automatically get back  the product
of the co-adoint orbits, or its diffeomorphic iamge.  Instead one gets the `twin' of the Cartesian product.
The concept `twin pair'  
was introduced in [C2, Sect. 12]. For a general compact $G$-symplectic
manifold, with the prequantum data $(M,L,f)$, such  that $f(M)$ is
transversal to $\ft$,  it has a `twin' $(M_G,L_G)$
which has  identical Riemann-Roch as that of $(M,L)$.  
In general $M_G$ is
only symplectic outside a set of real codimension 2. 

This raises the question  whether there is another construction of the fusion
product so that when the level $k\ge (\lda+\ldaa|\theta)$, one gets back 
the Cartesian product. Although we know a potential candidate for it, there
is some technical difficulty in proving the necessary properties.

\subsection{Relation with others work}
There have been a long list of work related to fusion product and 
 the moduli of flat
connections,   see [Be, MS, L, JK]. 

The present work has one  advantage that we have found a 
a compact holomorphic model ($X_N$) which has all the relavant information about $\cMab$.
And the fixed point formula proved in [C2] enables us to calculate 
 explicitly and directly the Riemann-Roch of $Y=G\times _T X_N$. The variety  $X_N$ serves naturally
as the compactified quotient of $\cMab$ by the nilpotent subgroup
 $LG^{\C+}$. In general  quotient of a variety by nilpotent group  rarely
exists, even in finite dimension.   
The existing work on moduli of  parabolic bundles, and on extended moduli spaces
of flat connections do not provide the fusion product or $X_N$.
 Hence the fusion tensor can not be realized  geometrically as it is done here.

Recent work [L]  is related in terms of calculating the Riemann-Roch
number of moduli of parabolic bundles, at least for the $SU(n)$-case. 
The calculation of other invariants, e.g symplectic volumes, was  worked
out there for general $G$. In that regard, [JK] is also related.
\subsection{ Acknowledgment} I am indebted to Prof.~Kac   who answered   my questions with great patience and insight, his suggestion to study the
$SU(3)$-case
 in detail was valuable.
   Prof.~Liu and I had many enjoyable and informative  conversations, I thank
him for that.

\include{2}

\include{3}

\include{4}

\include{5.final}

\include{bib}

\end{document}

%% file: 2.tex
\section{Moduli space of flat connections on a Riemann surface with boundary}
Let $S$ be a Riemann surface with boundary, the number of boundary components
is $n$. 
Let  $G$  be a connected, simply-connected and compact   simple Lie group.
Let $\cA, \cAp$ be the space of all smooth $G$-connections over  $S, \pS$ of a fixed
$G$-principle bundle respectively,  
and $\cG,\cGp$ be the group of smooth gauge transformations,  $F_A=d_AA$  the curvature
operator.  Because $S$ is of dimension 2, the
principle bundle can be assumed to be the trivial one. 
Clearly  $\cGp=LG^n$ where $LG$ is the loop group.

\subsection{What topology?} 
Let $|\cdot|_{l,S} , |\cdot|_{k,\pS}$   denote the norms of the Sobolev spaces $W^l(S), W^k(\pS)$ over $S$ and $\partial S$ respectively.

Let $\cA, \cG$   be completed by $|\cdot|_{l-1,S}+|\cdot|_{k-1,\pS}$
and $|\cdot|_{l,S}+|\cdot|_{k,\pS}$ respectively, denote the completions
by $\cAlk, \cGlk$.  

Then the gauge transformation is a bounded operator 
$$\cAlk\times  \cGlk \rightarrow \cAlk.$$ 

Let $\cApk, \cGpk$ be the completions of $\cAp, \cGp$ under $W^{k-1}(\partial
S),W^{k}(\partial S)$ norms respectively.

We make the observation that the restriction
\begin{equation}\begin{split}
&A\in \cAlk \mapsto A|_\pS\in  \cApk,\\
&g\in \cGlk \mapsto g|_\pS  \in \cGpk
\end{split}\end{equation}
 are  bounded operators.  Each element $d+ad\theta\in \cApk$  has an
extension to $S$  inside  $W^{k+1/2}(S)$, in particular it is in $\cAlk$ if
$l\leq k+1/2$. In
fact the harmonic  extension has this property by interior regularity result.
In other words, the restriction map defined above is onto. Likewise for $\cGpk$.

From now on, we assume that $k\geq 2, k+1/2\geq l>2$. 
This way, all the elements are continuous.

Applying the  proof of Theorem 3.2 of [FU] to both $S, \pS$, 
we conclude $\{g_i\}$  has subsequential convergence in $\cGlk$ if $g_i(A)\rightarrow B$ in $\cAlk$.
Hence the orbit space is Hausdorff.

The tangent space of $\cAlk$ is given by
\begin{equation}
T_A\cAlk=\{a\in  W^l(S, \fg\times T^*(S)
 \big| \,   a|_\pS \in \cApk\}.
\end{equation}

Let $\cGIlk$ be the subgroup of $\cGlk$ with $I$ on the $\pS$. 
The following can be verified in the same way as in the case $\pS=\emptyset$, utilizing the
condition that $\Lie \cGIlk$ has 0 boundary values on $\pS$:
\begin{lemma}
The curvature operator  $F: \cAlk \cA^{l-1,k-1}$ with $F(A)=F_A=d_AA$ is a moment map for the action
by $\cGlk$, with respect to the symplectic form $\Omega$ on $T\cAlk$
 given by
$$\Omega (a,b)=\frac{1}{2\pi}\int_S(a,b)$$
where  $(,)$ is the same  bilinear  form on $\fg$ as before.

If $*$ is the Hodge operator defined by the Riemann surface,  
then $\Omega(\cdot,*\cdot)$ is positive definite.
\end{lemma}

\begin{proposition}
1). The space of flat connections mod out the action by 
$\cGIlk$, $\cM=F^{-1}(0)/\cGIlk$,  
 is an infinite dimensional complex Hilbert manifold.
The tangent space at $A$ is the space   
$$T_{[A]}\cM=\{a\in T_A\cAlk|\, d_Aa=d_A*a=0\}.$$

2). The action by $(LG)^n\simeq \cGlk/\cGIlk\simeq \cGpk$ on $\cM$, defined by extending each $g\in (LG)^n$ to an element in $\cGlk$, $\tg$, then $g[A]=[\tg A]$, 
is Hamiltonian with respect to the induced two form on $T_{[A]}\cM$ from
$\Omega$. 
The moment map is given by $\Phi([A])=a$ if $A_\pS=d+ad\theta$ with $a\in 
(\lfg)^n$.

3). There is a line bundle $\cL$ over $\cM$  on which $\cLG$ acts covering the action
by $LG$ on $\cM$. It has a invariant connection whose curvature is $2\pi
i\Omega$. 
\end{proposition}

\subsection{Further symplectic reduction on $\cM$}
Let $0<m\le n$,    fix $m$ components of $\pS$, and the corresponding
$m$-components $\Psi$ in the moment map   $\Phi$. 
The action by different components of $LG$ in $(LG)^n$ commutes, as 
it can be observed from the definition of its action on $\cM$. 

\begin{figure}\label{21}
\begin{center}
\input Fig/1_1.tex
\caption{ m=3, n=5} 
\end{center}
\end{figure}

Clearly $\Psi$ is the moment map associated with the action by
$(LG)^m$ on the selected $m$-components in $\pS$. 
Let $\vecv\in (\lfg)^m$ be a fixed elements.
And $(LG)^m_\vecv$ be the stabilizer of $d+\vecv d\theta$.
\begin{proposition}
The space $\cMv=\Psi^{-1}(\vecv)/(LG)_\vecv^m$ is a complex Hilbert manifold 
on which $(LG)^{n-m}$ acts in a Hamiltonian manner. 
The line bundle $\cL$ descends to one on $\cM_v$ denoted by $L$ on which
$(\cLG)^{n-m}$ acts. The two form $\Omega$ descends to $\omega$ on $\cMv$.
It is the $2\pi i$-curvature of an $(\cLG)^{n-m}$-invariant connection. 
\end{proposition}

{\it Remark:} One can combine the two above to get $\cMv$ directly from $\cA$.
The action with respect to the full gauge transformation is Hamiltonian,
the moment map is given by $(F_A,\Phi(A))$. 
 Selecting $m$-boundary components. Then the action by the subgroup
$\cGlk_m$
 whose elements are    $I$ on the other $n-m$ components, is Hamiltonian with respect to
the map $(F,\Psi)$. 

Given $\vecv$,  define the subgroup $\cGvlk\subset \cGlk_m$  which consists
of elements with boundary values in $(LG)_\vecv^m$ on the $m$-designated
boundary components, and $I$ on the rest of boundary.
The new subgroup is the stabilizer of $(0, \vecv)$.
The quotient $(F^{-1}(0), \Psi^{-1}(\vecv))$ by $\cGvlk$ is the 
 desired space.

The proof of the above is a standard exercise.

Let's state an obvious fact:
\begin{lemma}
The map $\mu:\cMv\rightarrow (\lfg)^{n-m}$, $\mu([A])$ given by the boundary
values at  the other $n-m$ components is proper.
It is the moment map with respect to the action by $LG^{n-m}$.
\end{lemma}
This is obtained by noticing that $\Psi^{-1}(C^m)$ is compact where $C$ is 
the affine alcove of $\fg$.
\subsection{Transversality}

\begin{lemma}
Suppose  $0<m<  n$, then for generic value $\vecv\in C^m$,
$\cMv$ satisfies the transversality condition.
\end{lemma}
{\it Pf:} 
Write $\Psi=(\psi, \mu)\in \lfg^m\times \lfg^{n-m}$. 
First we show that the slice $\mu^{-1}(\ft^{n-m})$ is a smooth submanifold in $\cM$. We  verify the claim by showing if $\mu(A)\in (\ft)^{n-m}$, $D\mu(T_A\cM)$ is onto $(\lfg/\ft)^{n-m}$. 
Let $\fk=(\lfg)^n_\Psi$, the stabilizer of $\Psi(A)$. It is a direct product
of the stabilizer of each component in $\Psi(A)$. Hence it is finite dimensional.
So $Z_\fk=\Psi^{-1}(\fk) $ is  finite dimensional, smooth and symplectic near $A$.
And $d\Psi$ is onto  $(\lfg)^n/\fk$. 
Write $\fk=(\fk_1,\fk_2)\subset \lfg^m\times \lfg^{n-m}$.
We need to check that $D_A\mu(T_AZ_\fk)$ is onto  $\fk_2/(\fg)^{n-m}$.
Otherwise, from finite dimensional theory, 
there is a $\eta\in \fk_2/(\ft)^{n-m}$  fixing $A$, but if we take the
extension of $\eta$ as having 0 boundary value on the first $m$-components, 
we know that no flat connections having such a stabilizer.
Thus $D_A\mu(T_AZ_\fk)$ is onto  $\fk_2/(\fg)^{n-m}$.
Or $D_A\mu(T_A\cM)$ is  onto $\lfg/\ft$, and $\mu^{-1}(\ft^{n-m})$ is a smooth
submanifold in $\cM$. 
Apply Sard's theorem to $$\psi:\mu^{-1}(\ft^{n-m})\rightarrow \lfg^m,$$
thus for generic $\vecv$, $\psi^{-1}(\vecv)$ is smooth, which  also
implies  the action by $(LG)_\vecv^m$ on the first $m$-components
has only discrete stabilizer.  Therefore, we have the generic smoothness
of $\cMv$.

 We may further assume that $\vecv\in C^m$, otherwise
we may transform $\vecv$ there by conjugation, and the discussion above
about the smoothness of $\psi^{-1}(\vecv)$ remains valid   under 
  conjugation. QED

%% file: Fig/1_1.tex
\setlength{\unitlength}{0.00066667in}
\begingroup\makeatletter\ifx\SetFigFont\undefined%
\gdef\SetFigFont#1#2#3#4#5{%
  \reset@font\fontsize{#1}{#2pt}%
  \fontfamily{#3}\fontseries{#4}\fontshape{#5}%
  \selectfont}%
\fi\endgroup%
{\renewcommand{\dashlinestretch}{60}
\begin{picture}(2852,3645)(0,-10)
\drawline(610,336)	(645.512,401.812)
	(698.819,470.250)
	(763.966,527.812)
	(835.000,561.000)

\drawline(835,561)	(900.337,562.359)
	(974.327,545.625)
	(1060.000,486.000)

\drawline(1060,486)	(1007.090,415.453)
	(950.247,379.037)
	(881.406,346.875)
	(806.765,322.729)
	(732.520,310.359)
	(664.866,313.529)
	(610.000,336.000)

\drawline(610,336)	(610.000,336.000)

\drawline(1060,486)	(1032.641,556.541)
	(1010.860,618.141)
	(994.436,671.889)
	(983.151,718.875)
	(975.117,796.922)
	(985.000,861.000)

\drawline(985,861)	(1002.961,909.082)
	(1029.002,961.594)
	(1061.875,1014.809)
	(1100.331,1065.000)
	(1188.995,1141.406)
	(1236.704,1160.168)
	(1285.000,1161.000)

\drawline(1285,1161)	(1338.902,1129.131)
	(1374.536,1065.609)
	(1397.354,981.439)
	(1405.661,935.050)
	(1412.807,887.625)
	(1426.347,795.170)
	(1443.425,715.078)
	(1510.000,636.000)

\drawline(1510,636)	(1549.719,675.844)
	(1570.233,768.750)
	(1581.596,821.215)
	(1598.130,870.281)
	(1660.000,936.000)

\drawline(1660,936)	(1716.110,949.283)
	(1776.505,949.828)
	(1837.721,939.428)
	(1896.299,919.875)
	(1948.777,892.963)
	(1991.694,860.484)
	(2035.000,786.000)

\drawline(2035,786)	(1949.590,693.000)
	(1873.616,658.312)
	(1810.000,636.000)

\drawline(1810,636)	(1764.312,627.000)
	(1703.770,624.000)
	(1621.343,627.000)
	(1569.726,630.750)
	(1510.000,636.000)

\drawline(2035,861)	(1979.074,932.525)
	(1934.649,994.828)
	(1901.285,1048.998)
	(1878.543,1096.125)
	(1863.165,1173.609)
	(1885.000,1236.000)

\drawline(1885,1236)	(1942.630,1275.375)
	(2024.785,1283.250)
	(2112.047,1267.500)
	(2185.000,1236.000)

\drawline(2185,1236)	(2232.753,1167.188)
	(2246.960,1118.930)
	(2258.511,1068.000)
	(2285.014,977.812)
	(2335.000,936.000)

\drawline(2335,936)	(2376.439,965.953)
	(2402.121,1037.625)
	(2431.743,1114.734)
	(2485.000,1161.000)

\drawline(2485,1161)	(2531.199,1168.131)
	(2579.429,1166.109)
	(2672.950,1138.125)
	(2747.496,1084.078)
	(2785.000,1011.000)

\drawline(2785,1011)	(2784.424,959.062)
	(2768.050,905.625)
	(2738.502,853.500)
	(2698.405,805.500)
	(2650.382,764.438)
	(2597.057,733.125)
	(2541.055,714.375)
	(2485.000,711.000)

\drawline(2485,711)	(2399.403,765.750)
	(2366.310,833.906)
	(2350.652,880.395)
	(2335.000,936.000)

\drawline(610,336)	(568.280,413.909)
	(528.317,489.248)
	(490.082,562.088)
	(453.550,632.502)
	(418.691,700.561)
	(385.479,766.337)
	(353.887,829.902)
	(323.886,891.328)
	(295.450,950.688)
	(268.550,1008.052)
	(243.160,1063.493)
	(219.252,1117.084)
	(196.799,1168.895)
	(175.772,1219.000)
	(156.145,1267.469)
	(137.890,1314.375)
	(105.386,1403.785)
	(78.041,1487.807)
	(55.635,1567.014)
	(37.949,1641.984)
	(24.762,1713.292)
	(15.854,1781.514)
	(11.007,1847.224)
	(10.000,1911.000)

\drawline(10,1911)	(14.406,1982.832)
	(25.393,2060.080)
	(42.361,2141.807)
	(64.712,2227.078)
	(91.848,2314.956)
	(123.171,2404.506)
	(158.081,2494.791)
	(195.981,2584.875)
	(236.273,2673.823)
	(278.357,2760.697)
	(321.635,2844.563)
	(365.510,2924.484)
	(409.382,2999.525)
	(452.653,3068.748)
	(494.725,3131.219)
	(535.000,3186.000)

\drawline(535,3186)	(590.253,3249.844)
	(664.720,3320.250)
	(749.327,3379.781)
	(835.000,3411.000)

\drawline(835,3411)	(884.403,3413.742)
	(940.533,3410.625)
	(999.666,3401.508)
	(1058.076,3386.250)
	(1112.040,3364.711)
	(1157.831,3336.750)
	(1210.000,3261.000)

\drawline(1210,3261)	(1210.814,3212.344)
	(1195.584,3162.375)
	(1166.858,3113.344)
	(1127.185,3067.500)
	(1079.114,3027.094)
	(1025.193,2994.375)
	(967.973,2971.594)
	(910.000,2961.000)

\drawline(910,2961)	(837.593,2968.031)
	(768.145,3004.500)
	(694.624,3075.469)
	(654.140,3125.473)
	(610.000,3186.000)

\drawline(1210,3261)	(1198.269,3199.264)
	(1188.600,3141.785)
	(1181.018,3088.343)
	(1175.553,3038.719)
	(1171.079,2950.043)
	(1175.399,2874.000)
	(1188.731,2808.832)
	(1211.295,2752.781)
	(1285.000,2661.000)

\drawline(1285,2661)	(1355.146,2621.484)
	(1435.405,2610.375)
	(1515.462,2624.578)
	(1585.000,2661.000)

\drawline(1585,2661)	(1627.218,2713.172)
	(1648.035,2782.125)
	(1654.797,2861.391)
	(1654.855,2944.500)
	(1655.556,3024.984)
	(1664.248,3096.375)
	(1735.000,3186.000)

\drawline(1735,3186)	(1810.817,3166.312)
	(1879.056,3087.000)
	(1950.267,3000.938)
	(2035.000,2961.000)

\drawline(2035,2961)	(2084.805,2966.525)
	(2137.181,2981.953)
	(2189.085,3005.561)
	(2237.474,3035.625)
	(2311.529,3108.234)
	(2335.000,3186.000)

\drawline(2335,3186)	(2286.828,3256.266)
	(2240.905,3283.646)
	(2185.994,3305.625)
	(2126.271,3322.025)
	(2065.913,3332.672)
	(2009.097,3337.389)
	(1960.000,3336.000)

\drawline(1960,3336)	(1866.749,3298.125)
	(1808.323,3253.172)
	(1735.000,3186.000)

\drawline(2335,3186)	(2303.562,3128.531)
	(2273.568,3072.899)
	(2245.001,3019.051)
	(2217.843,2966.930)
	(2192.077,2916.482)
	(2167.687,2867.651)
	(2122.963,2774.625)
	(2083.534,2687.411)
	(2049.263,2605.570)
	(2020.012,2528.663)
	(1995.644,2456.250)
	(1976.021,2387.892)
	(1961.006,2323.148)
	(1950.463,2261.581)
	(1944.252,2202.750)
	(1942.238,2146.216)
	(1944.283,2091.539)
	(1950.249,2038.280)
	(1960.000,1986.000)

\drawline(1960,1986)	(1973.025,1936.023)
	(1989.778,1886.748)
	(2010.562,1837.845)
	(2035.678,1788.984)
	(2065.430,1739.837)
	(2100.118,1690.072)
	(2140.045,1639.362)
	(2185.514,1587.375)
	(2236.826,1533.783)
	(2294.283,1478.256)
	(2358.188,1420.464)
	(2428.842,1360.078)
	(2506.549,1296.768)
	(2591.609,1230.205)
	(2684.325,1160.059)
	(2733.649,1123.539)
	(2785.000,1086.000)

\dottedline{50}(1885,261)(1135,336)
\drawline(1257.390,353.911)(1135.000,336.000)(1251.419,294.208)
\dottedline{50}(2335,261)(2035,561)
\drawline(2141.066,497.360)(2035.000,561.000)(2098.640,454.934)
\dottedline{50}(2710,261)(2710,561)
\drawline(2740.000,441.000)(2710.000,561.000)(2680.000,441.000)
\drawline(1135,1986)	(1222.036,1967.953)
	(1297.638,1955.062)
	(1363.127,1947.328)
	(1419.820,1944.750)
	(1512.092,1955.062)
	(1585.000,1986.000)

\drawline(1585,1986)	(1652.389,2043.440)
	(1708.836,2130.555)
	(1734.605,2188.537)
	(1759.615,2257.893)
	(1784.527,2339.941)
	(1797.152,2386.137)
	(1810.000,2436.000)

\drawline(1285,1986)	(1357.901,2078.493)
	(1416.186,2143.309)
	(1465.128,2185.720)
	(1510.000,2211.000)

\drawline(1510,2211)	(1563.835,2219.708)
	(1604.868,2217.531)
	(1660.000,2211.000)

\drawline(535,1611)	(607.880,1561.712)
	(671.133,1519.789)
	(725.857,1484.572)
	(773.151,1455.401)
	(849.844,1412.563)
	(910.000,1386.000)

\drawline(910,1386)	(953.740,1371.175)
	(1006.801,1355.557)
	(1066.308,1340.321)
	(1129.390,1326.641)
	(1193.173,1315.694)
	(1254.784,1308.655)
	(1311.351,1306.698)
	(1360.000,1311.000)

\drawline(1360,1311)	(1451.871,1337.901)
	(1503.806,1361.179)
	(1561.791,1392.180)
	(1627.365,1431.784)
	(1702.066,1480.869)
	(1787.432,1540.315)
	(1834.594,1574.198)
	(1885.000,1611.000)

\drawline(985,1461)	(1055.895,1498.981)
	(1117.730,1530.588)
	(1171.603,1556.261)
	(1218.614,1576.440)
	(1296.441,1602.072)
	(1360.000,1611.000)

\drawline(1360,1611)	(1444.071,1595.621)
	(1504.792,1572.220)
	(1585.000,1536.000)

\drawline(385,2436)	(406.929,2373.095)
	(427.070,2319.206)
	(463.750,2234.962)
	(498.555,2176.238)
	(535.000,2136.000)

\drawline(535,2136)	(581.828,2106.345)
	(642.385,2082.101)
	(705.499,2066.057)
	(760.000,2061.000)

\drawline(760,2061)	(823.559,2069.928)
	(901.386,2095.560)
	(948.397,2115.739)
	(1002.270,2141.412)
	(1064.105,2173.019)
	(1135.000,2211.000)

\drawline(535,2286)	(617.700,2295.793)
	(679.251,2299.057)
	(760.000,2286.000)

\drawline(760,2286)	(822.543,2240.872)
	(861.399,2197.662)
	(910.000,2136.000)

\put(835,186){\makebox(0,0)[lb]{\smash{{{\SetFigFont{6}{7.0}{\rmdefault}{\mddefault}{\updefault}$\pS_1$}}}}}
\put(1660,486){\makebox(0,0)[lb]{\smash{{{\SetFigFont{6}{7.0}{\rmdefault}{\mddefault}{\updefault}$\pS_2$}}}}}
\put(2560,561){\makebox(0,0)[lb]{\smash{{{\SetFigFont{6}{7.0}{\rmdefault}{\mddefault}{\updefault}$\pS_3$}}}}}
\put(610,3486){\makebox(0,0)[lb]{\smash{{{\SetFigFont{6}{7.0}{\rmdefault}{\mddefault}{\updefault}$\pS_4$}}}}}
\put(1885,3411){\makebox(0,0)[lb]{\smash{{{\SetFigFont{6}{7.0}{\rmdefault}{\mddefault}{\updefault}$\pS_5$}}}}}
\put(1610,36){\makebox(0,0)[lb]{\smash{{{\SetFigFont{9}{10.0}{\rmdefault}{\mddefault}{\updefault}$\vecv=(v_1,v_2,v_3)$}}}}}
\end{picture}
}

%% file: 3.tex
\section{When $S$ is a pair of pants}
The situation  when $S$ is the sphere after removig three disjoint
disks is of special interest.

Let $m=2$, and $\vecv=(a,b)$ be a pair of elements in $\ft\times\ft$. 
From the last section, we know that for generic $\vecv$, $\cMv$ is a
comlex $LG$-manifold satisfying the transversality condition.

From [C1,C2] we have constructed a pair of  $T$ and $G$-orbifold $X_N, Y=G\times_T X_N$ respectively. So what are those spaces in the present situation?

Let $(x,y)=(\exp(\lp a),\exp(\lp b))$, $G^\sss_c$ denote  the semi-simple subgroup
commuting with $c\in T$. Naturaly $G^\sss_c\cap T$ is a  maximal torus. 

\begin{proposition}
 $X_N$ is  constructed from a subset of the Cartisian product of the conjugacy classes 
passing $a,b$ respectively,  then collapsing certain orbits of the subgroups in
$T$: 
\begin{equation}
X_N=\{ (r,s)\in \Ad_G x\times \Ad_G y|\, rs \in e^{2 \pi i C}\}/\simeq,
\end{equation}
where $\simeq$ is defined by $(r,s)\simeq (tr,ts)$ for $t\in
G^\sss_{c}\cap T$ where $c=rs$.

$Y=G\times_T X_N$   and it also  can be described as a  blow up of the following:
$$\{ (r,s)\in \Ad_G x\times \Ad_Gy\}/\simeq$$
where $\simeq$ is define by
$(r,s)\simeq (gr,gs)$ for $ g\in G^\sss_c	\cap \Ad_{c}T$ where $c=rs$.
\end{proposition}

With this description, one can see that $Y$ is related  but may not be
diffeomorphic 
to $\Ad_Ga\times \Ad_Gb$.

\begin{definition}
The fusion product of the  co-adjoint orbits $\Ad_Ga\times \Ad_Gb$ is defined to
be $Y$ as in the above.
\end{definition}
\subsection{Equivariant Riemann-Roch of  $Y$}
Suppose $\lda,\ldaa \in P_k^+$.  Suppose the pair $ (a,b)=(\lda/k,\ldaa/k)$
satifies the generic condition, 
the moduli space $\cMab$ and the line bundle $L^k$ induce the space 
$Y$ and a holomorphic $G$-line  bundle  $V$ over $Y$.

Let $\chi(g)$  be the character of the representation of $G$ on 
$H^0(Y,V^k)$, 
 the  equivariant Riemann-Roch of the pair be $\RR(Y,V^k)$.

\subsection{Fixed point sets on the fusion product}
Given a $\tau\in\explife$, let us identify the fixed points on $\cMab$.
Each $[A]\in \cMab$ represents a class of flat connection.  Choose
a  base point $p$  on the third boundary component and evaluate the holonomy of
the connection $A$ with that base point. One gets a pair $(x,y)\in
\Ad_Ge^{\lp a}\times\Ad_ Ge^{\lp b}$. 
Suppose  the image $\mu(A)$ is in $\ft$.
That $A$ is fixed by $\tau$  implies  $(x,y)$ is fixed by $\tau$ under
conjugation.
 Because $\tau$ is a regular element in $T$,  we know that
$x,y\in T$. Hence $(x,y)\in W(e^{\lp a})\times W(e^{\lp b})$. 
Therefore all the flat connections fixed by $\tau$ produces such
representations.  In this case, the $\tau$-fixed
points and $T$-fixed points are the same.

Suppose $[A']$  has the same holonomy representation.
Then there is  a gauge transformation  $g$ on $S$ such that
$g(p)\in T$, and $g(A')=A$. 
Since  $\mu([A]),\mu([A'])$ are fixed by $\tau$, they are in $\ft$. Or 
$$\mu([A])=\mu(g[A'])=\Ad_g\mu([A'])\in\ft.$$ One concludes that
$$g^{-1}dg+g^{-1}\mu g\in\ft,\quad g(p)\in T.$$ 
Therefore $g^{-1}(p)$ is an element in the co-root lattice.
Thus if we requie $\mu(A)\in W(C\setminus C^\aff)$, then there is only one $[A]$ which
gives the desired representation. The holonomy of $A$ around the 
the third boundary  component is given by $z$ with 
\begin{equation}
\label{31}
z=e^{\lp \mu(A)}=u( e^{\lp a})\cdot
v(e^{\lp b}).
\end{equation}

One knows further that the image  $\mu(A)\in W(C^\inte)$.
Otherwise, $\mu(A)$ has nontrivial $(LG )_\mu^\sss$, and clearly 
$T\subset (LG )_\mu$. We know $T\subset (LG)_A$, because 
$A$ is fixed by $\tau$ iff $A$ is fixed by $T$.
Therefore
$$[\lfg_\mu,\lfg_\mu]\cap \lfg_A\neq0,$$
which violates the transversality condition.

Now we end the discussion by stating
\begin{lemma}
The $\tau$-fixed points  on $\cMab$ are 1-1 correspondence with
the set $W(e^{\lp a})\times W(e^{\lp b})$. The image of the moment map of 
the element $[A]$ with $\mu(A)\in W(C^\inte)$ is the unique solution to
Eq.~\ref{31}
\end{lemma}

Such fixed points on $\cMab$ induces the same on $Y$. 
 In the current situation  there is no need to consider
the closure of the $\tau$-fixed points on $Y$ intersecting 
the compactfying locus, since all of them having images in $W(C^\inte)$. Therefore the remainder term  $\RT(\tau)$ defined in [C2] is 0  for
all $\tau\in \explife$.

%% file: 4.tex
\section{Weights at the fixed points and Riemann-Roch}

Compare the result on fixed points from the previous section, with the fixed points in the product of the coadjoint 
orbits, $\Ad_G a\times \Ad_Gb $, one sees the 1-1 correspondence between the 
fixed point sets  of $\Ad_G a\times \Ad_Gb $ and the interior 
fixed points of the fusion product.

For the product of  the coadjoint orbits, the image of a  fixed point  is 
$u(a)+v(b)$, while for the fusion product, the image is    $c\in W(C^\inte)$ with
$e^{\lp c}=u(e^{\lp a})\cdot v(e^{\lp b})$ or $c=u(a)+v(b)+t$ where $t$ is  
a translation element in $W^\aff$, with $e^{2\pi i t}=I$.  

In the following $p$ denote a  fixed point on the fusion product and $q$ the
corresponding fixed point on
the Cartesian product.

Let $A\in \cMab$ be a fixed point, we want to calculate the weights on
its tangent space. As the recipe provided in [C2] shows that the information
on the weights and the fixed points will determine the equivariant
Riemann-Roch of $Y$, hence the multiplicity of $\cLG$,  $G$-irreducible highest weight module in $H^0(\cMab, L^k), H^0(Y,V)$ respectively.

First let's examine the $G=SU(2)$ case. 
\begin{lemma}
1). For $G=SU(2)$, $C=[0,1 ]$. For $a,b\in C$,   the four fixed points on
$\cMab$ with images in $[-1,1]=W(C)$ are $$(e^{2\pi  ai},e^{2\pi bi}) , \, (e^{-2\pi a i},e^{2\pi bi}), \, (e^{2\pi a i},e^{-2\pi bi}) ,\,
(e^{-2\pi a i},e^{-2\pi bi}).$$ 
2). The weight on $T_AX_N$ which has $\dim_\C=1$,   is always given by $\pm \alpha$ where $\alpha$ is the root.
Two out of four fixed  points have  images in $C$.
Suppose $a\ge b$, then $c(e^{2\pi a i},e^{-2\pi bi})=a-b\in C $
the weight there is   $\alpha$.
  \begin{equation}
 \begin{cases}
 a+b< 1 :& \quad  c(e^{2\pi  ai},e^{2\pi bi})=  a+b\in C; \\
 a+b>    1: &
\quad  c(e^{-2\pi a i},e^{-2\pi bi})=2 -( a+b)\in  C,
\end{cases}
\end{equation}
the weights on $T_AX_N$ in both cases are given by $\alpha$.

\end{lemma}

{\it Pf:}  Part 1 is obvious.
\begin{figure}\label{41}
\begin{center}
\input Fig/3_1_new.tex
\caption{Weights   at the fixed points of Cartesian and fusion products for
$SU(2)$-case. The middle one is part of $\mu(\cMab)\cap \ft$ and the
bottom one is  $\phi(\Ad_Ga\times \Ad_Gb)\cap\ft$.}
\end{center}
\end{figure}

 2).       Fig.~4.1 illustrates what is going on  here.

{\it Remark:} The weights in the figure have  the opposite signs  from  what appear
below in the text, because the weights discussed here refer to the
weights by $t\in T$ on the tangent space, while what appear in the figure  are 
those   by the action of $t^{-}$. 

  The end points of
$C$,  $\alpha=0, 1$ are  identified with $0, \alpha^\vv/2$.    Since $e^0=I, e^{\pi i\alpha^\vv}=-I$ are in the center of the group,
the image of $\cMab$  does not intersect $\partial C$.
 Otherwise we have $rs$ is in the center of the group, 
and $r,s$ after conjugating by the same element in $G$, we may assume
$r,s\in T$. Therefore $A$ is fixed by $T$,
and we have $$(\lfg)_A\cap[\lfg_\mu,\lfg_mu]\supset \ft\cap su(2)=\ft$$
which violates the transversality condition.
Thus we have $\mu(A)\cap \partial C=\emptyset$.

Each end point of an  interval in  $\mu(A)\cap C^\inte$ must be the image of a fixed
point, from a basic  property of the  moment map. Only two fixed  points
have their images in $C^\inte$, as we have identified them. Thus the image is
a single  interval. (In fact the image $\mu(X)\cap C$ is always a convex polytope, an
observation I first made  in loop group setting. But we can prove
this without much work in the present situation.)

The local property near $\mu(A)$ determine  the signs of the
weights, which is a general fact from symplectic geometry. In fact, $\mu(A)$ is a left end point then the weight is  $-\alpha$, and it is a right end point, 
then the weight is $\alpha$. 

In both case of part 2), $a+b> a-b$, and $2 -(a+b)>a-b$. Therefore the
weight is always $\alpha$. QED

As shown by the above, the sign of the weight has something to do with the
length of $a+b$, it is important.

Next we write the above result  in a form more convient for later use. 

Let $u,v\in W$, and $t$ be a translation by an element in the   lattice
$M$  generated by $ W(\theta^\vv)$.  

\begin{corollary}
Let $\phi:\Ad_G a+\Ad_Gb\rightarrow \fg^*$ be the moment map of the Cartesian
product, and  $\phi(q)=u(a)+v(b)$ .
Suppose $ c(p)=u(a)+v(b)+t\in C$, then the weight  on $T_A X_N$ differs  
with that at $T_q  \phi^{-1}(\ft)$ iff $| <\alpha^\vv,u(a)+v(b)>|>1$. 
\end{corollary}
{\it Pf:} By assumption, we have $|<\alpha^\vv,u(a)>|, |<\alpha^\vv,v(b)>|<1$, and
$$0<\, <\alpha^\vv,c(p)>=<\alpha^\vv,u(a)+v(b)+t><1.$$
 We have seen that if
 $ <\alpha^\vv,u(a)>, <\alpha^\vv,v(b)>$ have different signs,
 then $|<\alpha^\vv,u(a)+v(b)>|<1$, and   the weight at $p$ is  $-\alpha$. On the other hand
  the weight of $T_q  \phi^{-1}(\ft)$ 
 is given by $-\alpha$ as well, as it can be checked
easily. Thus they have the same sign.

If $<\alpha^\vv,u(a)>, <\alpha^\vv,v(b)>$ have the same sign $+$, then 
from the prevous lemma, 
 the sign of the weight of $T_p\mu^{-1}(\ft)$  agree with that of $T_q  \phi^{-1}(\ft)\cap(
\Ad_G a+\Ad_Gb)$. 

If  both $<\alpha^\vv,u(a)>, <\alpha^\vv,v(b)>$ have sign $ -$, then
$\alpha(t)=2$, and $$|<\alpha^\vv,u(a)>+ <\alpha^\vv,v(b)>|>1.$$ From the lemma, the  weight on  $T_p\mu^{-1}(\ft)$ is $\alpha$ while 
at $T_q \phi^{-1}(\ft)$ it is $-\alpha$.

Hence the only time the weights on $T_p\mu^{-1}(\ft), T_q\phi^{-1}(\ft)$ are
different, 
for $\mu(p)=u(a)+v(b)+t$ and $\phi(q)=u(a)+v(b)$ is when the 
the said condition is met. QED

From Fig. 4.1, one can see in the middle figure, the intervals
$[-2+a+b,b-a], [a+b,2+b-a]$ looks the same. 
This is because inside $\cMab$, 
 $$\mu^{-1}([-2+a+b,b-a])=R_{\alpha^\vv}\mu^{-1}([a+b,2+b-a]),$$
where $R_{\alpha^\vv}$ defined by the translation element $\alpha^\vv\in
W^\aff$ acts on $\cMab$. 
It preserves the complex structure and commutes with 
the $T$-action, therefore the two subvarieties have the same weights at 
the corresponding points. And one can furthre tell the difference 
between the Cartesian and fusion products at the point with image $a+b$.

\subsection{General $G$}
Now we want to compare the weights at for $X_N$, or $Y=G\times _T X_N$ with
 the Cartesian product.
 Fig.~\ref{42} illustrates the difference in $\mu(Y)\cap \ft$ and
$\phi(\Ad_Ga\times \Ad_Gb)\cap \ft$. 
\begin{figure}\label{42}
\begin{center}
\input Fig/4_2.tex
\caption{ The  left one is  $\phi(\Ad_Ga\times \Ad_Gb)\cap \ft$, the
other is $\mu(Y)\cap \ft$.} 
\end{center}
\end{figure}

Here we further assume that $a, b\in C^\inte$.
Without that, e.g. $b$ is a vertax $\neq 0$, then $e^{2\pi ib}$ is in the 
center of the group. Thus $\Ad_Ge^{2\pi ib}=e^{2\pi ib}$, while
$\Ad_Gb$ has positive dimension. Therefore  it is hard to compare the two varieties when one of $a,b$ is on $\partial C$.

Obviously the condition is generic, and we will  remove it in future.
 
The previous corollary can now  be used to  deal with the general cases.

From [K,  Chap. 6], we know that $M$ as a group of translation and the
Weyl group of  $\fg$, $ W$,  generate  the affine Weyl group $W^\aff$.
And $M$ is a normal subgroup. Thus each element of $\ft$ can be translated by an
element in $M$  to $W(C)$.

Suppose $(e^{2\pi   u(a)i},e^{2\pi v(b)i})$ defines a $T$-fixed point $p$  with image in $C$, i.e. $\mu(p)=u(a)+v(b)+t\in C$.
Then $\mu(p)\in C^\inte$ by transversality.
Let $$q=(u(a),v(b))\in \Ad_Gu(a)\times \Ad_Gv(b).$$  
For each positive root $\alpha$, the subalgebra $\fg_\alpha$ induces 
the subvariety 
$\Ad_{G_\alpha}e^{2\pi i u(a)}\times \Ad_{G_\alpha} e^{\lp v(b)}$ 
in $\cMab$, its intersection with $\mu^{-1}(\ft)$ is symplectic at
$p$, since $\mu(p)\in C^\inte$.

  For $a,b\in C^\inte$, there is no semi-simple stabilizer of either $u(a),v(b)$.
Thus the dimension of the subvariety is four and its intersection with $\mu^{-1}
(\ft)$ is of dimension 2 whose tangent space has a weight $\pm \alpha$. 

The restriction of the two form $\Omega$ is positive definite to that tangent
subspace, because the tangent vectors there are $\fg_\alpha$-valued 1-forms,
and both $T$-action on the tangent space and the $*$-operator preserves the subspace at $p$.

The subvariety $\Ad_{G_\alpha}e^{2\pi i u(a)}\times
\Ad_{G_\alpha} e^{\lp v(b)}$ is a $G_\alpha$-variety. The restriction
of the moment map $\mu$ to the subvariety the moment map for the
$G_\alpha$-action.
Let $\ft_\alpha=\fg_\alpha\cap\ft$.  
Then for each pair $(r,s)\in\Ad_{G_\alpha}e^{2\pi i u(a)}\times
\Ad_{G_\alpha} e^{\lp v(b)}$, if $e^{\lp \mu}=rs$,  $\mu$ is on
the affine line $$\ft_\alpha+\mu(p)=\ft_\alpha+(u(a)+v(b)+t).$$
 Because the images pass $\mu(p)$  and must be in the direction of $\ft_\alpha$.
 As assumed $\mu(p)\in C$, thus $0<(\alpha|\mu(p))<1$.

On the other hand, inside the Cartesian product of the coadjoint orbits, we
have a similar subvariety $\Ad_{G_\alpha} u(a)\times \Ad_{G_\alpha}
v(b)$.  Its intersection with $\phi^{-1}(\ft)$ is also symlectic, and with
2-d  tangent space. The weight is given by $\pm\alpha$.

\begin{proposition}\label{key}
The weights on $$T_p(\Ad_{G_\alpha}e^{2\pi i u(a)}\times
\Ad_{G_\alpha} e^{\lp v(b)}),\quad  T_q(\Ad_{G_\alpha}u(a)\times
\Ad_{G_\alpha}  v(b))$$
differ iff $$|<\alpha^\vv, u(a)+v(b)>|>1.$$
\end{proposition}
{\it Pf:} 
Let $R_t$ denote the action by the translation element $t$ in $W^\aff$.
As commented after Cor.~4.1, inside $\cMab$,  
the weights at $ T_p$ and $T_{R_{-t}p}$ are the same. 
Obviously $R_{-t}p$ is a $T$-fixed point.

   The connection $A$ which defines $R_{-t}p$, has boundary value 
at the 3rd component given by $ \mu(A)=u(a)+v(b)$, since
$$\mu(R_{-t}p)=R_{-t}\mu(p)=R_{-t}(u(a)+v(b)+t)=u(a)+v(b)$$ by the equivariance of the moment map.  
Choose $a',b'\in \R_+ \simeq (\ft_{\alpha})_+$, so that 
 $$<\alpha^\vv,u(a)>=\pm a', \quad <\alpha^\vv,v(b)>=\pm b'.$$
Now $<\alpha^\vv,u(a)+v(b)> >1$ iff $a',b'$ have the same sign
and $|a'+b'|>1$.
Then  the variety $\Ad_{G_\alpha}e^{\lp a'}\times \Ad_{G_\alpha} e^{\lp b'}$
is exactly the same situation as the $SU(2)$-case, thus the weight has a
different sign from that of  $T_q$ iff $|a'+b'|>1$, or iff
$<\alpha^\vv,u(a)+v(b)> >1$.
QED
\subsection{Comparision of $\FC_p,\FC_q$ for general $G$}

 Next we find another way of writing the condition when the sign 
of the weights differ.

All the facts about affine Lie algebra $\fg^\aff$ based on $\fg$  can be found in [K, Chap.~6].
The Lie algebra $\fg^\aff$ is the  extention of  $\clg$ by differentiation,
$d$ on the circle.  
 The Lie algebra $\fg^\aff$    has $\ft\oplus \R K \oplus\R d$ as its
Cartan subalgebra. Here $K$ is the central element and $d$ is the
differentiation.  The dual is given by $\ft^*\oplus \R \Lambda_0\oplus\R
\delta$. All the real positive  affine  roots are of the form
$$ \Delta_+(\fg^\aff)=\Delta_+(\fg)\cup  \{ n\delta\pm\alpha |n\in \Z_+, \alpha\in
\Delta_+(\fg)\}.$$
 It turns out the condition for when the signs of the weights differ is 
best written  in the language of affine Lie algebras.

Assume  $c=u(a)+v(b)+t\in C$, where $t$ is in the coroot lattice,
hence  $0<\alpha(c)<1$, for all positive  roots.
As discussed earlier, the signs differ iff $|\alpha(u(a)+v(b))|>1$.

1).  $ \alpha(c-t)=\alpha(u(a)+v(b))<-1$.
It can be written as $ (\delta +\alpha)(R_{-t}c)<0$.
Thus the sign of $\alpha$ as a weight changes iff $(\delta +\alpha)(R_{-t}c)<0$ 
or equivalently $R_t(\delta +\alpha)<0$.
In this case, at $p$ the weight is $\alpha$, at $q$ it is $-\alpha$.

2). $ \alpha(c-t)=\alpha(u(a)+v(b))>1$.
It can be written as $ (\delta -\alpha)(R_{-t}c)<0$.  
In this case, the sign changes iff $R_{t}(\delta -\alpha)<0$. And the weight
is $-\alpha$ at $p$ and $\alpha$ at $q$.

Are there other affine roots changing sign under $R_t$? 
The anwer is no.
Because when  $n\geq 2$,
 $$n\pm<\alpha^\vv,u(a)+v(b)>\ge n-2\ge 0.$$
Therefore the sign does not change under $R_t$.

Thus we obtain the following:
\begin{lemma}
 1). The  sign changes   iff $R_t(\gamma)<0$ for a positive affine root
$\gamma=\delta\pm\alpha$.

2). 
Let the set of those roots whose induced weights on the  tangent spaces at
$p,q $ differ by sign be denoted by $S$.
Then $\sum_S\alpha=\sum_S \gamma \mod \delta=\sum_{R_t\gamma<0} \gamma \mod \delta$.
\end{lemma}

Next we compare $\det_{T_pY}(1-t^{-1})$ with  $\det_{T_q (\Ad_G a\times \Ad_G b)}(1-t^{-1})$. 

For $c  \in C$,  we have
\begin{equation}\label{denor}
\begin{split}
&\det_{T_pY}(1-t^{-1})=\prod_{\alpha\in \Delta_+(\fg)}(1-e^{-\alpha})\det_{T_p
\mu^{-1}(\ft) }(1-t^{-1});\\
&\det_{T_q(\Ad_G a\times \Ad_G b)}(1-t^{-1})=\prod_{\alpha\in \Delta_+(\fg)}(1-e^{-w\alpha})\det_{T_q
\phi^{-1}(\ft) }(1-t^{-1}),
\end{split}
\end{equation}
where  $w\in W$  in the above satisfies $\phi(q)=u(a)+v(b)=R_t(c)\in w(\ft_+)$. 

Comparing the two in Eq.~(\ref{denor}), we obtain for $s=e^{2\pi i\tau}$,
\begin{equation}
\begin{split}
 &\det_{T_q(\Ad_G a\times \Ad_G b)}(1-s^{-1})\\
&=(-1)^{|w|}e^{-\sum_{w\alpha<0}w\alpha}\prod_{\alpha\in \Delta_+(\fg)}(1-e^{-\alpha})
\det_{T_q \phi^{-1}(\ft) }(1-s^{-1})\\
&=(-1)^{|w|+\#S}e^{-\sum_{w\alpha<0}w\alpha}\prod_{\alpha\in \Delta_+(\fg)}(1-e^{-\alpha})
e^{\sum_{\alpha\in S} \alpha}\det_{T_p \mu^{-1}(\ft) }(1-s^{-1})\\
&=(-1)^{|w|+\#S}e^{-\sum_{w\alpha<0}w\alpha+\sum_{\alpha\in S}\alpha}
\det_{T_pY}(1-s^{-1}).
\end{split}
\end{equation}

We already know that 
\begin{equation}
\begin{split}
&\sum_{\alpha\in S} \alpha=\sum_{ R_t(\gamma)<0}\gamma\mod \delta;\\
&-\sum_{w\alpha<0,\alpha\in\Delta_+(\fg)}w\alpha= \sum_{w^{-1}\beta<0,\beta\in \Delta_+(\fg)}\beta.
\end{split}
\end{equation}
Now $w^{-1}\beta<0$ iff $R_t(\beta)<0$. To see that, write $\beta=-w\alpha$,
then $w^{-1}\beta<0$ iff $\alpha>0$;
 we also  know that  $$(R_t(\beta)|c)=(\beta|c-t)=-(w\alpha|c-t)$$
where $c-t\in w(\ft_+^*)$. Hence $R_t(\beta)<0$ iff $-(w\alpha|c-t)<0$, or iff
$-(\alpha|\ft_+^*)<0.$ Thus
$$-\sum_{w\alpha<0,\alpha\in \Delta_+(\fg)}w\alpha=\sum_{R_t\beta<0,\beta\in
\Delta_+(\fg)}\beta.$$
 Furthermore, we have checked
that no other positve affine roots change sign under $R_t$. 
Apply the well known formula:
$$(\rho+h^\vv\Lambda_0)- w^{-1}(\rho+h^\vv\Lambda_0)=\sum_{\gamma\in\Delta(\fg^\aff),w(\gamma)<0}\gamma,\quad \forall w\in W^\aff
$$
we obtain the following  key identity:
\begin{proposition}
\begin{equation}\label{affine}
\begin{split}
-\sum_{w\alpha<0}w\alpha+\sum_{\alpha\in S}\alpha&=\sum_{\beta\in
\Delta_+(\fg^\aff), R_t\beta<0}\beta \quad \mod \delta\\
&= -R^{-1}_t(\rho+h^\vv\Lambda_0)+(\rho+h^\vv \Lambda_0) \mod \delta\\
&= h^\vv t \mod \delta;\\
(-1)^{|w|+\#S}&= (-1)^{|R_t|}\\
&=1,
\end{split}
\end{equation}
where the absolute value sign denotes the length of an element in $W^\aff$.
( Translation elements  in $W^\aff$ have even  lengths).
\end{proposition}

Now it is easy to prove the main result:
\begin{theorem}
Let $\lda,\ldaa\in P_+^k$, hence $a=\lda/k, b=\ldaa/k\in C$.

1). Let $p,q$ denote the fixed points on $Y$,  $\Ad_G \lda\times \Ad_G\ldaa$
respectively, such that $c=\mu(p)=u(a)+v(b)+t\in C$, and $\phi(q)=u(a)+v(b)$.
Let $V$ be the  same line bundle as before on $Y$, and the $L_\lda\otimes L_\ldaa$ on
$\Ad_Ga\times \Ad_Gb$.
Denote by  $\FC_p, \FC_q$ respectively the contributions of $p,q$ to
$\RR(Y,V^k), \RR(\Ad_G \lda\times \Ad_G\ldaa,L_\lda\otimes L_\ldaa )$. Then we
have 
$$FC_p(e^{2\pi i\tau})=e^{2\pi i<(k+h^\vv)t,\tau>} \FC_q(e^{2\pi i\tau}).$$

2). $$\RR(Y,V^k)(s)=\chi_\lda(s)\cdot\chi_\ldaa (s)\onexplatk.$$
\end{theorem}
{\it Remark:} This  proves Verlinde conjecture.

3). Write $\RR(Y,V)(s)=\sum_{l\in P^k_+}m_l s^l$, then
$m_c$ is given by 
$$\frac{(-1)^l}{\big|\frac{M^*}{(k+h^\vv)M}\big|}\sum_{s\in \explife}
 \chi_{\bar{c}}(s)D^2(s)\chi_\lda(s)\chi_\ldaa(s). $$

{\it Pf:}  1).  The two functions $\FC_p, \FC_q$ can be written down as
\begin{equation}
\begin{split}
\FC_p(e^{ 2\pi i\tau})&=\frac{e^{ 2\pi ik< u(a)+v(b)+t,\tau>}}{\det_{T_pY}(1-e^{-2\pi i\tau})};\\
\FC_q(e^{ 2\pi i\tau})&=\frac{e^{ 2\pi ik < u(a)+v(b),\tau> }}
{\det_{T_q (\Ad_G a\times \Ad_G b)}(1-e^{-2\pi i\tau})}
\end{split}
\end{equation}
we already have compared the denominators in  the above, using Eq.~(\ref{affine})
to  obtain 
\begin{equation}
\begin{split}
\FC_p&=e^{2\pi i h^\vv <t,\tau>}\frac{e^{ 2\pi ik< u(a)+v(b)+t,\tau>}}{\det_{T_q (\Ad_G a\times \Ad_G b)}(1-t^{-2\pi i\tau})}\\
&=e^{2\pi i <(h^\vv+k)t,\tau> }\FC_q,
\end{split}
\end{equation}
the two sides are equal when $<(h^\vv+k) t,\tau>\in \Z$, or when $$t\in\latk.$$
 Thus we have the conclusion.

2). Using fixed points formula proved in [C2],
we obtain 
\begin{equation}
\begin{split}
\RR(Y,V)(e^{ 2\pi i\tau})&= \sum_{p} \FC_p(e^{ 2\pi i\tau})\onlatk,\\
\end{split}
\end{equation}
which is amazing considering  that $Y$ has lots more fixed points 
than those interior ones. See figure on the right  in fig. 4.2. 
All the fixed points with images on the boundary of $wC$ of which there are plenty, do not appear in the above. 
Thus  we have by Part 1),
\begin{equation}
\begin{split}
\RR(Y,V)(e^{ 2\pi i\tau})&= \sum_{p} \FC_p(e^{ 2\pi i\tau})\onlatk,\\
&=\sum_q \FC_q(e^{ 2\pi i\tau}) \onlatk\\
&=\RR(\Ad_G a\times \Ad_G b, L^k_a\otimes L^k_b )\onlatk\\
&=\RR(\Ad_G \lda\times \Ad_G\ldaa, L_\lda\otimes L_\ldaa )\onlatk\\
&=\chi_\lda\cdot\chi_\ldaa \onlatk.
\end{split}
\end{equation}

3). The expression follows right away from the expression of $\RR(Y,V^k)$ and
Cor. 1.1 of [C2]. 
QED

%% file: Fig/3_1_new.tex
\setlength{\unitlength}{0.00055555in}
\begingroup\makeatletter\ifx\SetFigFont\undefined%
\gdef\SetFigFont#1#2#3#4#5{%
  \reset@font\fontsize{#1}{#2pt}%
  \fontfamily{#3}\fontseries{#4}\fontshape{#5}%
  \selectfont}%
\fi\endgroup%
{\renewcommand{\dashlinestretch}{30}
\begin{picture}(7046,1920)(0,-10)
\drawline(3012,261)(1212,261)
\drawline(3612,261)(5412,261)
\drawline(5412,261)(4962,261)
\drawline(5082.000,291.000)(4962.000,261.000)(5082.000,231.000)
\drawline(3612,261)(4062,261)
\drawline(3942.000,231.000)(4062.000,261.000)(3942.000,291.000)
\drawline(1212,261)(1662,261)
\drawline(1542.000,231.000)(1662.000,261.000)(1542.000,291.000)
\drawline(3012,261)(2562,261)
\drawline(2682.000,291.000)(2562.000,261.000)(2682.000,231.000)
\drawline(12,1461)(6612,1461)
\drawline(6492.000,1431.000)(6612.000,1461.000)(6492.000,1491.000)
\drawline(3312,1611)(3312,1461)
\drawline(5112,1611)(5112,1461)
\drawline(1512,1611)(1512,1461)
\drawline(3012,861)(1812,861)
\drawline(3012,861)(2562,861)
\drawline(2682.000,891.000)(2562.000,861.000)(2682.000,831.000)
\drawline(1812,861)(2262,861)
\drawline(2142.000,831.000)(2262.000,861.000)(2142.000,891.000)
\drawline(3612,861)(4812,861)
\drawline(3612,861)(4062,861)
\drawline(3942.000,831.000)(4062.000,861.000)(3942.000,891.000)
\drawline(4812,861)(4362,861)
\drawline(4482.000,891.000)(4362.000,861.000)(4482.000,831.000)
\drawline(5412,861)(6612,861)
\drawline(5412,861)(5862,861)
\drawline(5742.000,831.000)(5862.000,861.000)(5742.000,891.000)
\drawline(6612,861)(6162,861)
\drawline(6282.000,891.000)(6162.000,861.000)(6282.000,831.000)
\put(3312,1761){\makebox(0,0)[lb]{\smash{{{\SetFigFont{5}{6.0}{\rmdefault}{\mddefault}{\updefault}$0$}}}}}
\put(1212,1761){\makebox(0,0)[lb]{\smash{{{\SetFigFont{5}{6.0}{\rmdefault}{\mddefault}{\updefault}$\alpha=-1$}}}}}
\put(4812,1761){\makebox(0,0)[lb]{\smash{{{\SetFigFont{5}{6.0}{\rmdefault}{\mddefault}{\updefault}$\alpha=1$}}}}}
\put(1062,411){\makebox(0,0)[lb]{\smash{{{\SetFigFont{5}{6.0}{\rmdefault}{\mddefault}{\updefault}$-a-b$}}}}}
\put(1212,36){\makebox(0,0)[lb]{\smash{{{\SetFigFont{5}{6.0}{\rmdefault}{\mddefault}{\updefault}$\alpha$}}}}}
\put(2637,36){\makebox(0,0)[lb]{\smash{{{\SetFigFont{5}{6.0}{\rmdefault}{\mddefault}{\updefault}$-\alpha$}}}}}
\put(2862,336){\makebox(0,0)[lb]{\smash{{{\SetFigFont{5}{6.0}{\rmdefault}{\mddefault}{\updefault}$b-a$}}}}}
\put(3462,336){\makebox(0,0)[lb]{\smash{{{\SetFigFont{5}{6.0}{\rmdefault}{\mddefault}{\updefault}$a-b$}}}}}
\put(5262,411){\makebox(0,0)[lb]{\smash{{{\SetFigFont{5}{6.0}{\rmdefault}{\mddefault}{\updefault}$a+b$}}}}}
\put(4962,36){\makebox(0,0)[lb]{\smash{{{\SetFigFont{5}{6.0}{\rmdefault}{\mddefault}{\updefault}$-\alpha$}}}}}
\put(3612,36){\makebox(0,0)[lb]{\smash{{{\SetFigFont{5}{6.0}{\rmdefault}{\mddefault}{\updefault}$\alpha$}}}}}
\put(2637,711){\makebox(0,0)[lb]{\smash{{{\SetFigFont{5}{6.0}{\rmdefault}{\mddefault}{\updefault}$-\alpha$}}}}}
\put(1662,711){\makebox(0,0)[lb]{\smash{{{\SetFigFont{5}{6.0}{\rmdefault}{\mddefault}{\updefault}$\alpha$}}}}}
\put(3612,711){\makebox(0,0)[lb]{\smash{{{\SetFigFont{5}{6.0}{\rmdefault}{\mddefault}{\updefault}$\alpha$}}}}}
\put(6462,711){\makebox(0,0)[lb]{\smash{{{\SetFigFont{5}{6.0}{\rmdefault}{\mddefault}{\updefault}$-\alpha$}}}}}
\put(5412,711){\makebox(0,0)[lb]{\smash{{{\SetFigFont{5}{6.0}{\rmdefault}{\mddefault}{\updefault}$\alpha$}}}}}
\put(4437,936){\makebox(0,0)[lb]{\smash{{{\SetFigFont{5}{6.0}{\rmdefault}{\mddefault}{\updefault}$2-a-b$}}}}}
\put(4437,711){\makebox(0,0)[lb]{\smash{{{\SetFigFont{5}{6.0}{\rmdefault}{\mddefault}{\updefault}$-\alpha$}}}}}
\put(1512,936){\makebox(0,0)[lb]{\smash{{{\SetFigFont{5}{6.0}{\rmdefault}{\mddefault}{\updefault}$-2+a+b$}}}}}
\put(2637,936){\makebox(0,0)[lb]{\smash{{{\SetFigFont{5}{6.0}{\rmdefault}{\mddefault}{\updefault}$b-a$}}}}}
\put(3612,936){\makebox(0,0)[lb]{\smash{{{\SetFigFont{5}{6.0}{\rmdefault}{\mddefault}{\updefault}$a-b$}}}}}
\put(5412,936){\makebox(0,0)[lb]{\smash{{{\SetFigFont{5}{6.0}{\rmdefault}{\mddefault}{\updefault}$a+b$}}}}}
\put(6237,936){\makebox(0,0)[lb]{\smash{{{\SetFigFont{5}{6.0}{\rmdefault}{\mddefault}{\updefault}$b-a+2$}}}}}
\end{picture}
}

%% file: Fig/4_2.tex
\setlength{\unitlength}{0.00052500in}
\begingroup\makeatletter\ifx\SetFigFont\undefined%
\gdef\SetFigFont#1#2#3#4#5{%
  \reset@font\fontsize{#1}{#2pt}%
  \fontfamily{#3}\fontseries{#4}\fontshape{#5}%
  \selectfont}%
\fi\endgroup%
{\renewcommand{\dashlinestretch}{30}
\begin{picture}(10165,4877)(0,-10)
\dashline{60.000}(2342,2431)(2331,4850)
\dashline{60.000}(2342,4831)(4442,3631)
\drawline(2512,4381)(3625,4381)
\drawline(3991,3732)(3625,4381)
\dottedline{75}(3617,4381)(2717,2806)
\drawline(2492,4381)(2342,4156)
\drawline(2192,4381)(2342,4156)
\drawline(2342,3406)(2717,2806)
\drawline(2717,2806)(3017,2806)
\drawline(3017,2806)(2867,2581)
\drawline(3992,3706)(3692,3181)
\drawline(3767,3181)(4292,3181)
\dottedline{75}(3242,3181)(3692,3181)
\dottedline{75}(3242,3181)(3017,2806)
\dottedline{75}(3692,3181)(3467,2806)
\dottedline{75}(3092,2806)(3467,2806)
\dottedline{75}(2342,4156)(2117,3706)
\dottedline{75}(2342,4156)(2567,3706)
\dottedline{75}(2117,3706)(2342,3406)
\dottedline{75}(2342,3406)(2567,3706)
\dottedline{75}(2567,3706)(3992,3706)
\dottedline{75}(2492,4381)(3242,3181)
\dashline{60.000}(5364,1231)(5353,3650)
\drawline(6930,2568)(6780,2828)
\drawline(5513,3222)(6117,3225)
\dottedline{75}(6342,2835)(6117,3225)
\dottedline{75}(6342,2835)(6780,2828)
\drawline(5678,1586)(5341,2170)
\drawline(5678,1586)(5948,1568)
\drawline(6598,1944)(6930,2568)
\dottedline{75}(5948,1568)(6225,1988)
\dottedline{75}(6598,1944)(6225,1988)
\dottedline{75}(6225,1988)(5513,3222)
\dottedline{75}(5678,1586)(6342,2835)
\dottedline{75}(5342,2131)(5567,2581)
\dottedline{75}(5567,2581)(5342,2956)
\dottedline{75}(5342,2956)(5492,3182)
\drawline(5492,3255)(5342,2956)
\dottedline{75}(5567,2581)(6917,2580)
\drawline(5342,2206)(5342,2956)
\drawline(5828,1326)(6165,742)
\dashline{60.000}(5364,1231)(7453,12)
\drawline(5828,1326)(5948,1568)
\drawline(7305,1919)(7455,1659)
\drawline(7163,365)(7467,887)
\dottedline{75}(7242,1276)(7467,887)
\dottedline{75}(7242,1276)(7455,1659)
\dottedline{75}(5948,1568)(6450,1599)
\dottedline{75}(6598,1944)(6450,1599)
\dottedline{75}(6450,1599)(7163,365)
\dottedline{75}(5828,1326)(7242,1276)
\dottedline{75}(6132,762)(6634,732)
\dottedline{75}(6634,732)(6847,350)
\dottedline{75}(6847,350)(7117,367)
\dottedline{75}(6634,732)(7309,1902)
\drawline(6197,725)(6847,350)
\dashline{60.000}(9520,1231)(9531,3650)
\drawline(7954,2568)(8104,2828)
\drawline(9371,3222)(8767,3225)
\dottedline{75}(8542,2835)(8767,3225)
\dottedline{75}(8542,2835)(8104,2828)
\drawline(9206,1586)(9543,2170)
\drawline(9206,1586)(8936,1568)
\drawline(8286,1944)(7954,2568)
\dottedline{75}(8936,1568)(8659,1988)
\dottedline{75}(8286,1944)(8659,1988)
\dottedline{75}(8659,1988)(9371,3222)
\dottedline{75}(9206,1586)(8542,2835)
\dottedline{75}(9542,2131)(9317,2581)
\dottedline{75}(9317,2581)(9542,2956)
\dottedline{75}(9542,2956)(9392,3182)
\drawline(9392,3255)(9542,2956)
\dottedline{75}(9317,2581)(7967,2580)
\drawline(9542,2206)(9542,2956)
\drawline(9056,1326)(8719,742)
\dashline{60.000}(9520,1231)(7431,12)
\drawline(9056,1326)(8936,1568)
\drawline(7579,1919)(7429,1659)
\drawline(7721,365)(7417,887)
\dottedline{75}(7642,1276)(7417,887)
\dottedline{75}(7642,1276)(7429,1659)
\dottedline{75}(8936,1568)(8434,1599)
\dottedline{75}(8286,1944)(8434,1599)
\dottedline{75}(8434,1599)(7721,365)
\dottedline{75}(9056,1326)(7642,1276)
\dottedline{75}(8752,762)(8250,732)
\dottedline{75}(8250,732)(8037,350)
\dottedline{75}(8037,350)(7767,367)
\dottedline{75}(8250,732)(7575,1902)
\drawline(8687,725)(8037,350)
\dashline{60.000}(7442,4831)(5342,3631)
\drawline(7067,2806)(6767,2806)
\drawline(5792,3706)(6092,3181)
\dottedline{75}(6542,3182)(6767,2806)
\drawline(7292,4381)(6617,4381)
\drawline(7292,4381)(7442,4156)
\drawline(7442,3405)(7067,2806)
\dottedline{75}(7442,4156)(7217,3706)
\dottedline{75}(7442,3405)(7217,3706)
\dottedline{75}(7217,3706)(5792,3706)
\dottedline{75}(7292,4381)(6542,3182)
\dottedline{75}(6652,4400)(6375,3980)
\dottedline{75}(6375,3980)(5937,3987)
\dottedline{75}(5937,3987)(5817,3744)
\drawline(5753,3708)(5937,3987)
\dottedline{75}(6375,3980)(7050,2811)
\drawline(6587,4362)(5937,3987)
\drawline(7592,4381)(8266,4381)
\dashline{60.000}(7442,4831)(9542,3631)
\drawline(7592,4381)(7442,4156)
\drawline(7442,3405)(7817,2806)
\drawline(7817,2806)(8117,2806)
\drawline(9092,3706)(8792,3181)
\dottedline{75}(8342,3182)(8792,3181)
\dottedline{75}(8342,3182)(8117,2806)
\dottedline{75}(7442,4156)(7667,3706)
\dottedline{75}(7442,3405)(7667,3706)
\dottedline{75}(7667,3706)(9092,3706)
\dottedline{75}(7592,4381)(8342,3182)
\dottedline{75}(8232,4400)(8509,3980)
\dottedline{75}(8509,3980)(8947,3987)
\dottedline{75}(8947,3987)(9067,3745)
\drawline(9131,3707)(8947,3987)
\dottedline{75}(8509,3980)(7834,2811)
\drawline(8297,4362)(8947,3987)
\dottedline{75}(6542,3182)(6092,3181)
\drawline(6092,3256)(6092,3256)
\drawline(6767,2806)(6167,3181)
\drawline(7442,4156)(7442,3406)
\drawline(8792,3181)(8117,2806)
\drawline(8267,1981)(8942,1531)
\drawline(7442,1681)(7442,856)
\dashline{60.000}(7442,4831)(7442,4231)
\dashline{60.000}(7442,3406)(7442,1681)
\dashline{60.000}(7442,856)(7442,31)
\dashline{60.000}(5417,3631)(6017,3256)
\dashline{60.000}(8942,1531)(9542,1231)
\dashline{60.000}(6842,2806)(8267,1981)
\dashline{60.000}(8117,2806)(6542,1906)
\dashline{60.000}(5342,1231)(5942,1531)
\dashline{60.000}(8867,3256)(9542,3631)
\drawline(6542,1906)(5942,1531)
\drawline(6617,1906)(7292,1906)
\drawline(7592,1906)(8267,1906)
\drawline(6842,331)(7142,331)
\drawline(7742,331)(8042,331)
\dashline{60.000}(2342,4831)(242,3631)
\dashline{60.000}(2342,2431)(2353,12)
\dashline{60.000}(2342,2431)(4442,1231)
\dashline{60.000}(2342,31)(4442,1231)
\dashline{60.000}(2337,2422)(248,3641)
\dashline{60.000}(2342,31)(242,1231)
\drawline(2172,4381)(1059,4381)
\drawline(693,3732)(1059,4381)
\dottedline{75}(1067,4381)(1967,2806)
\dashline{60.000}(242,3631)(242,1231)
\dashline{60.000}(4442,3631)(4442,1231)
\dashline{60.000}(4442,3631)(242,1231)
\drawline(568,1603)(12,2567)
\drawline(391,3209)(12,2567)
\dottedline{75}(16,2560)(1830,2568)
\drawline(4116,1603)(4672,2567)
\drawline(4293,3209)(4672,2567)
\dottedline{75}(4668,2560)(2854,2568)
\drawline(3946,1309)(3389,345)
\drawline(2644,352)(3389,345)
\dottedline{75}(3393,352)(2479,1919)
\drawline(738,1309)(1295,345)
\drawline(2040,352)(1295,345)
\dottedline{75}(1291,352)(2205,1919)
\drawline(2342,3406)(1967,2806)
\dottedline{75}(3317,1606)(3692,1606)
\drawline(3767,1606)(4142,1606)
\drawline(3767,1606)(3917,1306)
\dottedline{75}(3542,1981)(3767,1606)
\dottedline{75}(3542,1981)(3167,1981)
\drawline(3167,1981)(2492,1981)
\drawline(2867,2581)(3167,1981)
\dottedline{75}(3167,1981)(3317,1606)
\dottedline{75}(2342,1681)(2567,1306)
\dottedline{75}(2567,1306)(2342,856)
\drawline(2192,1906)(2342,1681)
\drawline(2342,856)(2042,331)
\drawline(2342,856)(2642,331)
\drawline(2492,1906)(2342,1681)
\dottedline{75}(2117,1306)(2342,1681)
\dottedline{75}(2117,1306)(2342,856)
\dottedline{75}(917,1606)(1292,1606)
\dottedline{75}(1292,1606)(1517,1981)
\dottedline{75}(1142,1981)(917,1606)
\dottedline{75}(1142,1981)(1517,1981)
\drawline(1517,1981)(2192,1981)
\drawline(1817,2506)(1517,1981)
\drawline(617,1606)(917,1606)
\drawline(917,1606)(767,1306)
\drawline(392,3181)(992,3181)
\drawline(692,3706)(992,3256)
\drawline(1667,2806)(1967,2806)
\drawline(1667,2806)(1817,2581)
\dottedline{75}(992,3181)(1517,3181)
\dottedline{75}(1292,2806)(1667,2806)
\dottedline{75}(1067,3181)(1292,2806)
\dottedline{75}(1667,2806)(1442,3181)
\dottedline{75}(2192,4381)(1442,3181)
\dottedline{75}(1292,2806)(617,1681)
\dottedline{75}(392,3181)(1142,1981)
\dottedline{75}(1367,1606)(2042,331)
\dottedline{75}(767,1306)(2117,1306)
\dottedline{75}(2567,1306)(3917,1306)
\dottedline{75}(3317,1606)(2642,331)
\dottedline{75}(3542,1981)(4292,3181)
\dottedline{75}(692,3706)(2117,3706)
\dottedline{75}(3467,2806)(4067,1606)
\drawline(8567,4006)(8942,4306)
\drawline(8641.963,4104.389)(8567.000,4006.000)(8679.445,4057.537)
\drawline(1292,331)(1517,781)
\drawline(1490.167,660.252)(1517.000,781.000)(1436.502,687.085)
\drawline(1292,331)(1742,331)
\drawline(1622.000,301.000)(1742.000,331.000)(1622.000,361.000)
\drawline(1292,331)(1067,781)
\drawline(1147.498,687.085)(1067.000,781.000)(1093.833,660.252)
\drawline(8492,4006)(8867,4006)
\drawline(8747.000,3976.000)(8867.000,4006.000)(8747.000,4036.000)
\drawline(8492,4006)(8267,4381)
\drawline(8354.464,4293.536)(8267.000,4381.000)(8303.015,4262.666)
\drawline(8492,4006)(8342,3706)
\drawline(8368.833,3826.748)(8342.000,3706.000)(8422.498,3799.915)
\put(917,181){\makebox(0,0)[lb]{\smash{{{\SetFigFont{5}{6.0}{\rmdefault}{\mddefault}{\updefault}$r_3(a)+r_3(b)$}}}}}
\put(8042,4581){\makebox(0,0)[lb]{\smash{{{\SetFigFont{5}{6.0}{\rmdefault}{\mddefault}{\updefault}$r_3(a)+r_3(b)+\alpha_3^\vv$}}}}}
\end{picture}
}

%% file: 5.final.tex
\section{ The induced representation on the fusion product and the proof
of an analogue of a conjecture by G. Segal  }
using the previous calculation of $\RR(Y,V^k)$. In order to do so, we need to
know something about the higher cohomology groups.
Ideally we want to show $H^i(Y,V^k)=0, i>0$. Or even better, the dual of the 
canonical line bundle $K^*(Y)$ is positive. That turns out to be more
involved. Instead, we will consider the reduced space of $Y$, or of $X_N$, and
use known results on the moduli of parabolic bundles.

\subsection{Comparison of complex structures}
Recall $\mu:X_N\rightarrow kC\subset \ft$ is a moment map with respect to a two 
form with degeneracy along $k\partial C$. 
For $a\in kC$, let $X_c=\mu^{-1}(c)/T$, and $L_c$ be  the corresponding 
line bundle over it. 
 Since $X_N$ is an orbifold, for generic value of $c$, 
$X_c$ is an orbifold.
From the characterization of $X_N$ in terms of the representation variety,
it is easy to see that 
$X_c, c\in kC^\inte$  is the same as that
of representation variety $$\{(h,k)\in \Ad_G e^{2\pi ia}\times \Ad_G e^{2\pi
ib}|hk=e^{2\pi ic/k}\}/T.$$
That variety has a complex structure via its well known diffeomorphism with   
the moduli of parabolic bundles.

\begin{proposition}
The two complex structure agree on $X_c$.
\end{proposition}
{\it Pf:}
First we describe  the complex structure on $X_c$ inherited from $X_N$.
Below we shall use $D\mu, D_\phi$ to denote the differential of $\mu,\phi$
so as not to confuse with the tangent vectors which are 1-forms.

What is $D_{(x,y)}(\mu-\muX)$? Recall that $T_ q\X=\lfg/\ft\oplus T_zX_\fg$
from
[C1], where $q=([I,z])$.
Let $\phi:X_\fg\rightarrow C\subset  \ft$ as in [C1], be the moment map of the
toric variety $X_\fg$. 

The tangent space to $T_{[p,q]}X_c$ is given by
$$T_{[p,q]}X_c=\{(x,y)\in T_{p}\cMab\times \lfg/\ft| D_{x}\mu =D_y\muX, D_{*x}\mu =D_{Jy}\muX, D_y\phi=D_{Jy}\phi=0\} .$$
As done in [C1], the complex structure $J$ on $ \lfg/\ft$ is given by
$Jy=-y^J$ where
$y^J$ is the standard complex structure on $\lfg/\ft$,
i.e. $y^J+iy$ is the boundary value of a  holomorphic function on the unit
disk.

We expand $y\in\lfg$ in terms of  $E_\gamma, E_{-\gamma}$.
The   positive affine roots $\gamma\in \Delta_+(\fg^\aff)$ are 
$$\{n\delta\pm\alpha| n\ge 1,\alpha\in \Delta_+(\fg)\}\cup \Delta_+(\fg^\aff)\cup\{n\delta|n\ge 1 \}$$
and accordingly $E_\gamma=z^n E'{\pm\alpha},  E'_\alpha$ or $z^n h_\alpha$,
where $E'_\alpha,h'_\alpha$ are Chevalley basis of $\fg$.

The  real form defining $\lfg$ has basis  given by $$x_\gamma=i(E_\gamma+E_{-\gamma}),\quad y_\gamma=E_\gamma-E_{-\gamma}.$$
By insisting on $y^J+iy$ being holomorphic, i.e. in $\sum _{\gamma>0}\C
E_\gamma$,  we obtain 
$$x_\gamma^J=-y_\gamma, \quad y_\gamma^J=x_\gamma.$$

Let $y=\sum a_\gamma x_\gamma+b_\gamma y_\gamma$, then
\begin{equation}\begin{split}
y^J&=\sum -a_\gamma y_\gamma+b_\gamma x_\gamma\\
&=-\sum (a_\gamma-ib_\gamma)E_\gamma-(a_\gamma+ib_\gamma)E_{-\gamma},
\end{split}\end{equation}
and 
\begin{equation}\label{51}
\begin{split}
D_{Jy}\mu&=-D_{y^J}\mu=-\sum \gamma(\mu)(a_\gamma-ib_\gamma)E_\gamma+(a_\gamma+ib_\gamma)\gamma(\mu)E_{-\gamma}\\
&=\sum a_\gamma |\gamma(\mu)| x_\gamma+b_\gamma |\gamma(\mu)|y_\gamma\\
\end{split}\end{equation}
where  the fact  $\gamma(\mu)$ is purely imaginary is used.

On the other hand,  the complex structure derived from the parabolic bundles can be described
by the following.
Let $S^1\times [0, \infty)$ be the conformal equivalent of the disk, with
coordinate $(\theta,u)$, and $*du=-d\theta, *d\theta=du$. 
A 1-form $a$ is in the tangent space to the parabolic bundle iff $
d_Aa=d_A*a=0$ on the extended surface $S\cup S^1\times [0, \infty)$, and
$a$ is in $L^2$.  

Given $x,y$ as in the above, we first extend $y$ as a harmonic section with respect
to $d_A*d_A$ on the cynlinder.
Let $y=\sum a_\gamma x_\gamma+b_\gamma y_\gamma$, and 
define $$\ty=\sum e^{-|\gamma(u)|}(a_\gamma x_\gamma+b_\gamma y_\gamma),$$ 
it is easy to check from this expression the extension is both harmonic (w.r.t
$d_A*d_A$)  and in $L^2$.  Now the coefficients of $d\theta$ in  $d_A\ty|_\pS$ agrees with 
$d_x\mu$ by assumption.  By direct calculation, the coefficent of $du$ in
$d_A\ty$ is given by
$$\sum a_\gamma |\gamma(\mu)| x_\gamma+b_\gamma |\gamma(\mu)|y_\gamma $$
which agree with the expression in Eq.~\ref{51},
hence  $x,d_A\ty$ agree on $\pS$ in both $d\theta, du$ directions. Therefore
they define a $L^2$-harmonic 1-form with respect to $d_A$  on the extended surface. 
From this, we conclude that the two complex structures agree.
QED

\subsection{Vanishing result}
The positivity of $K^*(X_c)$ implies
$H^i(X_c,L_c)=0, i>0$
for any semi-positive line bundle over $X_c$. We want  to  show 
\begin{equation}\label{52}
\tr(s|H^0(X_N,L_N))=\RR(X_N,L_N)
\end{equation}
 which holds if $H^i(X_N,L_N)=0, i>0$ or having a positive $K^*(X_N)$. 
 The proof of  the positivity of $K^*(X_N)$  turns out to be more
involved.  So we shall by pass that to show Eq.~\ref{52}.

 We have $\dim H^0(X_c,L_c)=\RR(X_c,L_c)$ which 
equals  the coefficient $m_c$ of the character of weight $c$ in $\RR(X_N,L_N)$, by the Abelian version of the result in [M]. 
We also know for $X_N$,  that 
$\dim H^0(X_c,L_c)$ equals the coefficient of the character of weight $c$ in
$\tr(s|H^0(X_N,L_N))$, as shown in [GS].  
Therefore, $$\tr(s|H^0(X_N,L_N))=\sum m_c s^c.$$
Apply the holomorphic induction to $Y=X\times_T X_N$, we know that
\begin{equation}\label{53}
\begin{split}
\tr(s|H^0(Y,L))&=\sum_W w \frac{\tr(s|H^0(X_N,L_N))}{\prod_{\alpha\in\Delta_+( \fg)}(1-s^{-\alpha})}\\
&=\sum_c m_c\chi_c.
\end{split}
\end{equation}
Now the Cor.~1.1 of [C2] gives an eplicit formula for $m_c$ in terms of
$\RR(Y,L)$:
\begin{equation}\label{54}
\begin{split}
m_c&=
\frac{(-1)^l}{\big|\frac{M^*}{(k+h^\vv)M}\big|}\sum_{s\in \explife}
 \chi_{\bar{c}}(s)D^2(s)\RR(Y,L)(s)\\
&=\frac{(-1)^l}{\big|\frac{M^*}{(k+h^\vv)M}\big|}\sum_{s\in \explife}
 \chi_{\bar{c}}(s)\chi_\lda(s)\cdot\chi_\ldaa(s)D^2(s).\\
\end{split}
\end{equation}
Thus the summation over $c\in P^k_+$  yields
\begin{equation}\label{55}
\begin{split}
\tr(s|H^0(Y,L))&=\sum_c m_{c\in P^k_+}\chi_c(s)\\
&=\frac{(-1)^l}{\big|\frac{M^*}{(k+h^\vv)M}\big|}\sum_{c\in P^k_+,s\in \explife}
\chi_c(s) \chi_{\bar{c}}(s)\chi_\lda(s)\cdot\chi_\ldaa(s)D^2(s)\\
&=\chi_\lda(s)\cdot\chi_\ldaa(s) \onexplife,
\end{split}
\end{equation}
where we have used the following fact from [K, Chap. 13]:
$$\frac{(-1)^l}{\big|\frac{M^*}{(k+h^\vv)M}\big|}\sum_{c\in  P^k_+}
m_c \chi_c(s) \chi_{\bar{c}}(s)D^2(s)=1.$$

Therefore we have  verified the following conjectured by Verlinde in [V]:
\begin{proposition}
$$\chi_{[\lda]\qtensor [\ldaa]}(s)=\tr(s|H^0(Y,L)) =\chi_\lda(s)\cdot\chi_\ldaa(s)
\onexplife.$$
\end{proposition}

For such $\lda,\ldaa$, after applying the result from [C1] and Eq. \ref{53}, 
we also obtain 
\begin{corollary}
 $$H^0(\cMab,L^k)=\sum_c m_c \tchi_{(c,k)}$$
\end{corollary}
which is  analogous to a conjecture by G. Segal, see [T].
The original conjecture is made for certain moduli of  holomorphic bundles instead of moduli of flat connections. Their equivalence can be established using
argumenent similar to  that of  Donaldson  for closed surface.

%% file: bib.tex
\section{References}
\noindent [AB1] M .\,F.\,Atiyah and R.\,Bott.{\it  A Lefschetz fixed
point formula for elliptic complexes. I} {\em Ann. Math.}  86
(1967), 374--407.

\noindent [AB2] M.\,F.\,Atiyah and R.\,Bott. {\it A Lefschetz fixed
point formula for elliptic complexes. II} Applications {\em Ann. Math.}
{\bf 87} (1968), 451--491.

\noindent [AB] M. Atiyah \& R. Bott, {\it The Yang-Mills equations over
Riemann
surfaces}, Phil. Trans. R. Soc. London, A308 (1982) 523.

\noindent[AS] M.\,F.\,Atiyah and G.\,B.\,Segal. {\it The index
of elliptic operators II}. {\em Ann. Math.} { 87} (1968), 531--545.

\noindent [AS1] M.\,F.\,Atiyah and I.\,M.\,Singer. {\it The
index of elliptic operators I.} {\em Ann. Math.} {\bf 87} (1968),
484--530.

\noindent [AS3]  M.\,F.\,Atiyah and I.\,M.\,Singer. {\it The
index of elliptic operators III.} {\em Ann. Math.} {\bf 87} (1968),
546--604.

\noindent [AS4] M.\,F.\,Atiyah and I.\,M.\,Singer. {\it The index of
elliptic operators IV.}  {\em Ann. Math.} {\bf 93} (1971), 119--138.

\noindent [BGV] N. Berline, E. Getzler \& M. Vergne, {\it Heat kernels and
Dirac operators}, preliminary version (1990).

\noindent [Be] A. Beauville, {\it Conformal blocks, fusion rules and the
Verlinde formula}. Preprint, (1994)

\noindent [Bi] J.-M.\,Bismut. {\it The infinitesimal Lefschetz
formulas: a heat equation proof}. {\em J. Funct. Anal.} {\bf 62} (1985),
435--457.

\noindent [B] R.\,Bott. {\it  Homogeneous vector bundles.} {\em Ann. Math.} {\bf
66} (1957), 203--248.

\noindent [C]  S.\,Chang. {\it A fixed point formula on orbifolds}.
MIT preprint 1995.

\noindent [C1] S.\,Chang. {\it Geometric  loop group actions and the highest weight modules}.  MIT preprint 1996.

\noindent [C2] S.\,Chang. {\it 
Fixed point formula and loop group action}.  MIT preprint 1997.

\noindent [F] G. Faltings, {\it A proof for the Verlinde formula}, J.
Algebraic
 Geom. 3, 347-374 (1994).

\noindent [JK] L. Jeffrey \& F. Kirwan, {\it Intersection pairings in moduli
spaces of vector bundles of arbitrary rank over a Riemann surface}
alg-geom/9608029.

\noindent [K] V. Kac.{\it Infinite dimensional Lie algebras} 3rd ed. Cambridge
University Press, 1990.

\noindent [KP] V. Kac \& D. Peterson, {\it Infinite-dimensional Lie algebras,
theta functions and modular forms,} Advances in Math. 53 (1984), 125-264

\noindent [KW] V. Kac \& M. Wakimoto, {\it Modular and conformal invariance
constraints in representation theory of affine algebras}, Advances in Math. 70
(1988), 156-234.

\noindent [KS] M. Kreck \& S. Stolz, {\it Nonconnected moduli spaces of positive
curvature metrics.} JAMS, Vol. 6, No. 4, 825-850 (1993).

\noindent [GS] V. Guillemin \& S. Sternberg, {\it Geometric quantization and
multiplicity of group representations}, Invent. Math. 67 (1982) 515-538.

\noindent [L] K. Liu, {\it Heat kernel and moduli spaces II} to appear in
Mathematical Research Letter.

\noindent [Ma] I. Macdonald. {\it  Affine Lie algebras and modular forms}, 
258-276, Seminaire
Bourbaki Vol. 1980/81 Exposes 561-578, Lect. Notes in Math. 901.

\noindent [M] E.\,Meinrenken. {\it Symplectic Surgery and the Spin-C Dirac
operator}, dg-ga 9504002. 

\noindent [MS] G. Moore \& N. Seiberg, {\it Classical and quantum conformal
field theory} Comm. Math. Phys. 123 (1989) 177-254.

\noindent [PS] A. Pressley \&  G. Segal, {\it Loop group}, Oxford Sciece
Publications, 1986.

\noindent [S] G. Segal, {\it Two dimensional conformal field theories and
modular functors}, IXth International Congress on Math. Phys. 1988, edited by B.
Simon
 {\it et al} 22-37.

\noindent [T] C. Teleman, {\it Lie algebra cohomology and the fusion rules.}
Comm. Math. Phys. 173 (1995), no. 2, 265--311.

\noindent [TZ] Y. Tian \& W. Zhang, {\it
Symplectic reduction and analytic localization }, Courant Institute Preprint (1996)

\noindent [TUY] A. Tsuchiya, K. Ueno \& Y. Yamada, {\it Conformal field theory
on universal family of stable curves with gauge symmetris}, Advanced Studies in
Pure Math. 19 (1982)

\noindent [Ve]  M. Vergne, {\it  Multiplicity formula for geometric
quantization
 Part I, II},
Duke Math. J. 82 (1996), 143-179; Duke Math. J. 82 (1996), 181-194.

\noindent [V] E. Verlinde, {\it Fusion rules and modular transformations in 2D
conformal field theory,} Nuclear Phys. B300 (1989) 360-376.

\vskip .3in
\noindent Department of Mathematics 2-271\\
 M. I. T.\\
Cambridge, MA 02139

\noindent schang@@math.mit.edu